\documentclass[a4paper,11pt,twoside,reqno]{amsart}

\usepackage[utf8]{inputenc}
\usepackage[plainpages=false,pdfpagelabels=true]{hyperref}
\usepackage{amssymb,amsthm,wasysym}
\usepackage{slashed}
\usepackage[margin=1in]{geometry}

\numberwithin{equation}{section}

\newtheorem{Satz}{Theorem}[section]
\newtheorem{Prop}[Satz]{Proposition}
\newtheorem{Lem}[Satz]{Lemma}

\newtheorem{Cor}[Satz]{Corollary}

\theoremstyle{definition}
\newtheorem{Dfn}[Satz]{Definition}
\newtheorem{Bem}[Satz]{Remark}

\newcommand{\tr}{\operatorname{Tr}}
\newcommand{\dv}{\text{dvol}_g}

\allowdisplaybreaks[1]

\renewcommand{\epsilon}{\varepsilon}

\newcommand{\N}{\ensuremath{\mathbb{N}}}
\newcommand{\R}{\ensuremath{\mathbb{R}}}

\newcommand{\sff}{\mathrm{I\!I}}

\title{Harmonic maps with torsion}
\author{Volker Branding}
\date{\today}
\address{University of Vienna, Faculty of Mathematics\\
Oskar-Morgenstern-Platz 1, 1090 Vienna, Austria\\}
\email{volker.branding@univie.ac.at}
\begin{document}

\begin{abstract}
In this article we introduce a natural extension of the well-studied equation for harmonic maps
between Riemannian manifolds by
assuming that the target manifold is equipped with a connection that is metric but has non-vanishing torsion.
Such connections have already been classified in the work of Cartan.
The maps under consideration do not arise as critical points of an energy functional 
leading to interesting mathematical challenges.
We will perform a first mathematical analysis of these maps which we will call harmonic maps with torsion.
\end{abstract} 

\maketitle

\section{Introduction and results}
The \emph{harmonic map equation} is one of the most studied partial differential equation for maps between Riemannian manifolds.
Given a smooth map \(\phi\) between two Riemannian manifolds \((M,g)\) and \((N,h)\) the harmonic map
equation can be obtained by computing the first variation of the energy of a map which is given by
\begin{align}
\label{energy}
E(\phi)=\int_M|d\phi|^2\dv.
\end{align}
The critical points of \eqref{energy} are characterized by the vanishing of the so-called \emph{tension field}
\begin{align}
\label{harmonic}
0=\tau(\phi):=\tr_g\nabla^{\phi^\ast TN}d\phi.
\end{align}
The (standard) harmonic map equation \eqref{harmonic} is a semilinear second order elliptic partial differential
equation for which many results on existence and qualitative behavior of its solutions could be achieved
over the years. For an overview on the current status of research we refer to the survey article \cite{MR2389639}.

In the literature on harmonic maps one usually chooses to utilize the Levi-Civita connection
on the target manifold \(N\).
If one defines harmonic maps via a variational principle, as we have done above,
then, as we will see later, the variational approach chooses to employ 
the Levi-Civita connection on the target manifold.

This article is devoted to a first study of harmonic maps that are coupled to a torsion endomorphism
on the target manifold. Harmonic maps with a connection different from the Levi-Civita connection
on the domain manifold have already been investigated in great generality.
Such kind of maps became known as \emph{V-harmonic maps} and include the classes of Hermitian, affine and Weyl harmonic maps into Riemannian manifolds,
see the introduction of \cite{Chen2019} for more details.

In order to obtain a generalization of the harmonic map equation that takes into account
a connection with metric torsion on the target manifold we have to abandon the variational
point of view. Although this leads to a straightforward generalization of the standard harmonic map equation \eqref{harmonic}
the non-variational nature of harmonic maps with torsion leads to severe technical difficulties.

The central equation studied in this article is 
\begin{align}
\label{harmonic-torsion}
\tau^{\scriptscriptstyle tor}(\phi):=\tau(\phi)+A(d\phi,d\phi)=0,
\end{align}
where the superscript ``tor'' represents that we are considering a connection with torsion.
The non-vanishing torsion is given by the torsion endomorphism \(A(d\phi,d\phi)\) on which we will
give more details later.

Solutions of \eqref{harmonic-torsion} will be called \emph{harmonic maps with torsion}.
Equation \eqref{harmonic-torsion} can be obtained by taking the standard harmonic map equation
and then changing to a connection with metric torsion.

The energy of a map \eqref{energy} is invariant under various symmetry operations
as for example the invariance under diffeomorphisms on the domain. By Noether's theorem these
invariances all lead to a conserved quantity which can then be successfully applied 
in the analysis of the harmonic map equation.

The heat flow method is a strong tool in the analysis of the standard harmonic map equation 
that allows to derive various existence results, see \cite{MR2431658} for an overview. 
Unfortunately, it seems that this method is no longer applicable (at least without making many additional assumptions)
in the presence of torsion. One should expect such a phenomena: The heat flow for standard harmonic maps
decreases the energy of a map \eqref{energy} and tries to flow it to a critical point.
As there does not exist an energy functional for harmonic maps with torsion the heat flow cannot find a direction
in which it decreases the energy and we cannot expect the gradient flow to converge.

While connections with metric torsion have been intensively studied in the physics literature,
so far, not many mathematicians have taken up this direction of research.
Let us mention several mathematical results that are connected to the equation
studied in this article.
\begin{enumerate}
\item \emph{Geodesics with vectorial torsion} were studied in \cite{MR2103403}.
Vectorial torsion comprises a particular kind of torsion endomorphism in \eqref{harmonic-torsion}.
\item 
 \emph{Dirac-harmonic maps with torsion}, which are a mathematical version of a model from supersymmetric quantum field theory, 
 have been investigated in \cite{MR3493217}.
 \item 
 \emph{VT-harmonic maps} \cite{Chen2019} are solutions of the following equation
\begin{align}
\label{vt-harmonic}
\tau(\phi)+d\phi(V)+\tr_g T(d\phi,d\phi)=0.
\end{align}
Here, \(V\in\Gamma(TM)\) and \(T\) is a (1,2)-tensor on \(TN\).
Although the form of \(VT\)-harmonic maps is very similar to harmonic maps with torsion
it seems that the additional structure in \eqref{harmonic-torsion} arising due to the torsion on \(N\)
leads to a rich mathematical framework.

\item 
\emph{Magnetic harmonic maps} \cite{MR3573990,MR3884770} from a two-dimensional domain share some similarities with \eqref{harmonic-torsion}. 
They arise as critical points of 
\begin{align}
\label{energy-magnetic}
E_B(\phi)=\int_M(|d\phi|^2+\frac{1}{2}\phi^\ast B)\dv,
\end{align}
where \(B\) is a two-form on \(N\) and \(\dim N\geq 3\).
The critical points of \eqref{energy-magnetic} are given by
\begin{align*}
\tau(\phi)=Z(d\phi(e_1)\wedge Zd\phi(e_2))
\end{align*}
and it is obvious that they have a structure analog to harmonic maps with torsion \eqref{harmonic-torsion}.
Here, \(\{e_1,e_2\}\) is an orthonormal basis of \(TM\) and 
the vector-bundle homomorphism \(Z \in \Gamma (\operatorname{Hom}(\Lambda^2T^\ast N,TN))\) is defined by the equation
\begin{align*}
\Omega(\eta,\xi_1,\xi_2)=\langle Z(\xi_1\wedge\xi_2),\eta\rangle,
\end{align*}
where \(\Omega=dB\) is a three-form on \(N\) and \(\xi_1,\xi_2,\eta\in\Gamma(TN)\).

\item In the study of \emph{pseudoharmonic maps} one is also naturally led to consider the case
of connections having torsion \cite{MR3938843,MR3844509}. However, these
maps still arise from a variational principle and the non-vanishing
torsion is given on the domain manifold.

\item The results of this article may also be interesting in the study of harmonic
maps to Lie groups as one can have a 
non-symmetric connection in the case of a non-abelian Lie group.
\end{enumerate}

In the case of \(N=\R^q\) with the flat metric the standard harmonic map equation reduces 
to the linear Laplace equation, whereas the equation for harmonic maps with torsion \eqref{harmonic-torsion}
would still be nonlinear due to the non-vanishing torsion.

Throughout this article the notation will be as follows:
Local coordinates on \(M\) will be denoted by \(x^i\) whereas
for local coordinates on \(N\) we will use \(y^\alpha\).
We will use Latin indices \(i,j,k=1,\ldots,m:=\dim M\) on the domain manifold \(M\)
while we will employ Greek letters \(\alpha,\beta,\gamma,\ldots,n:=\dim N\)
to represent indices on the target manifold \(N\).
Moreover, we will make use of the Einstein summation convention and tacitly
sum over repeated indices.

This article is organized as follows:

In section two we provide the necessary background on metric connections with torsion
and study the equation for harmonic maps with torsion in more detail.

The third section is devoted to geometric aspects of harmonic maps with torsion.
We derive various Bochner formulas, study the effects of conformal transformations
on harmonic maps with torsion and also discuss their stability.

In the last section we study several analytic aspects of harmonic maps with torsion.
We show that they satisfy the unique continuation property, study the regularity of weak solutions,
prove a removable singularity theorem and finally provide a Liouville type result
under a small energy assumption.

We will mostly apply techniques from standard harmonic maps 
in our analysis as far as they are still applicable.

\section{Harmonic maps with torsion}
Before we turn our attention to harmonic maps coupled to torsion let us give a short 
introduction to connection with metric torsion,
where we follow the introduction of \cite{MR2103403}.

\subsection{A shortcut to connections with metric torsion}
We consider a Riemannian manifold \((N,h)\) and by \(\nabla^{\scriptscriptstyle LC}\) we denote
its Levi-Civita connection. 
For any affine connection there exists a \((2,1)\)-tensor field \(A\) such that
\begin{align}
\label{nabla-torsion}
\nabla_XY=\nabla^{\scriptscriptstyle LC}_XY+A(X,Y)
\end{align}
for all vector fields \(X,Y\in\Gamma(TN)\).

We require that the connection \(\nabla\) is \emph{orthogonal}, that is for all vector fields \(X,Y,Z\) we have
\begin{align}
\label{nabla-metric}
\partial_X\langle Y,Z\rangle =\langle\nabla_XY,Z\rangle+\langle Y,\nabla_XZ\rangle,
\end{align}
where \(\langle\cdot,\cdot\rangle\) denotes the scalar product of the metric \(h\).
Combing \eqref{nabla-torsion} and \eqref{nabla-metric} it follows that the endomorphism \(A(X,\cdot)\) has to be skew-adjoint
\begin{align}
\label{torsion-skewadjointess}
\langle A(X,Y),Z\rangle=-\langle Y,A(X,Z)\rangle.
\end{align}

Any torsion tensor \(A\) induces a \((3,0)\)-tensor via the assignment
\[
A_{XYZ}=\langle A(X,Y),Z\rangle.
\]
We define the space of all admissible torsion tensors on \(T_pN\) by
\[
\mathcal{T}(T_pN):=\big\{A\in\otimes^3T^*_pN\mid A_{XYZ}=-A_{XZY}~ X,Y,Z\in T_pN\big\}.
\]
For \(A\in\mathcal{T}(T_pM)\) and \(Z\in T_pM\) we set
\[
c_{12}(A)(Z)=A_{\partial_{y^i}\partial_{y^i}Z},
\]
where \(\partial_{y^i}\) is a local basis of \(TN\).
Metric connections with torsion have been classified by Cartan \cite{MR1509255}
who proved the following
\begin{Satz}[Cartan, 1924]
\label{torsion-classification}
Assume that \(\dim N\geq 3\). Then the space \(\mathcal{T}(T_pN)\) has the following irreducible
decomposition
\[
\mathcal{T}(T_pN)=\mathcal{T}_1(T_pN)\oplus\mathcal{T}_2(T_pN)\oplus\mathcal{T}_3(T_pN),
\]
which is orthogonal with respect to \(\langle\cdot,\cdot\rangle\) and is explicitly given by
\begin{align*}
\mathcal{T}_1(T_pN)=& \{A\in\mathcal{T}(T_pN)\mid\exists V \textrm{ such that } A_{XYZ}=\langle X,Y\rangle\langle V,Z\rangle-\langle X,Z\rangle\langle V,Y\rangle\}, \\
\mathcal{T}_2(T_pN)=& \{A\in\mathcal{T}(T_pN)\mid A_{XYZ}=-A_{YXZ}~~\forall X,Y,Z\}, \\
\mathcal{T}_3(T_pN)=& \{A\in\mathcal{T}(T_pN)\mid A_{XYZ}+A_{YZX}+A_{ZXY}=0\textrm{ and } c_{12}(A)(Z)=0\}. 
\end{align*}
Moreover, for \(\dim N=2\) we have
\[
\mathcal{T}(T_pN)=\mathcal{T}_1(T_pN).
\]
\end{Satz}
For a proof of the above Theorem we refer to \cite[Theorem 3.1]{MR712664}.

\begin{Cor}
For any orthogonal connection on a Riemannian manifold \(N\) and \(\dim N>3\) the torsion tensor can be written as
\begin{align}
\label{torsion-decomposition}
A(X,Y)=\langle X,Y\rangle V-\langle V,Y\rangle X
+T(X,Y,\cdot)^\sharp+S(X,Y,\cdot)^\sharp,
\end{align}
where \(T(X,Y,\cdot)^\sharp\) and \(S(X,Y,\cdot)^\sharp\) are uniquely defined via
\[
T(X,Y,Z)=\langle T(X,Y,\cdot)^\sharp,Z\rangle, \qquad S(X,Y,Z)=\langle S(X,Y,\cdot)^\sharp,Z\rangle.
\]
\end{Cor}

We will use the following terminology:
A torsion endomorphism is called
\begin{enumerate}
\item \emph{vectorial} if it is contained in \(\mathcal{T}_1(T_pN)\),
\item \emph{totally antisymmetric} if it lies in \(\mathcal{T}_2(T_pN)=\Lambda^3T^\ast_pN\),
\item \emph{of Cartan type} if it is an element of \(\mathcal{T}_3(T_pN)\).
\end{enumerate}

Throughout this article we make use of the following sign convention for the curvature of a connection
\begin{align*}
R(X,Y)Z=\nabla_X\nabla_YZ-\nabla_Y\nabla_XZ-\nabla_{[X,Y]}Z
\end{align*}
for given vector fields \(X,Y,Z\).

We have the following relation between the curvature tensors of \(\nabla^{\scriptscriptstyle Tor}\) and 
\(\nabla^{\scriptscriptstyle LC}\)
\begin{align}
\label{curvature-tensors}
R^{\scriptscriptstyle Tor}(X,Y)Z=&R^{\scriptscriptstyle LC}(X,Y)Z
+(\nabla^{\scriptscriptstyle LC}_XA)(Y,Z)-(\nabla^{\scriptscriptstyle LC}_YA)(X,Z) \\
\nonumber &+A(X,A(Y,Z))-A(Y,A(X,Z)).
\end{align}	
It is obvious that the curvature tensor of a connection with metric torsion is antisymmetric in \(X,Y\)
as is the standard Riemann curvature tensor. However, not all of the symmetries of the Riemann 
curvature tensor still hold true in the presence of non-vanishing torsion.

The torsion tensor \(T(X,Y)\) is related to the torsion endomorphism \(A(X,Y)\)
via
\begin{align}
\label{def-torsion-tensor}
\nabla_XY-\nabla_YX-[X,Y]=T(X,Y)=A(X,Y)-A(Y,X)
\end{align}
for given vector fields \(X,Y\).

For more details on the geometry of Riemann manifolds having a connection with metric torsion we refer to the lecture notes \cite{MR2322400}.
Let us also mention the following geometric results connected to this article:
The uniformization theorem on closed surfaces for a metric connection with torsion
was proved via the Ricci flow in \cite{MR3667423}.
Geometric aspects of manifolds having a connection with vectorial torsion were studied in \cite{MR3457391}.

\subsection{Harmonic maps with torsion}
The way we have obtained the equation for harmonic maps with torsion \eqref{harmonic-torsion}
was by considering the standard harmonic map equation
and passing over to a connection with metric torsion on the target manifold \(N\).

We would like to emphasize once more that harmonic maps with torsion are non-variational and cannot be obtained as a critical point of an energy functional in general.

More precisely, if we consider a variation of the map \(\phi\) given by \(\phi_t\colon M\times (-\epsilon,\epsilon)\to N\) satisfying \(\frac{\partial\phi_t}{\partial t}\big|_{t=0}=\eta\),
where \(\epsilon>0\) is a small number,
a direct calculation shows
\begin{align*}
\frac{d}{dt}\big|_{t=0}\frac{1}{2}\int_M|d\phi_t|^2\dv
=&\int_M\langle\frac{\nabla^{\scriptscriptstyle Tor}}{dt}d\phi_t(e_i),d\phi_t(e_i)\rangle\dv\big|_{t=0} \\
=&\int_M\big(\langle\frac{\nabla^{\scriptscriptstyle LC}}{dt}d\phi_t(e_i),d\phi_t(e_i)\rangle
+\underbrace{\langle A(d\phi_t(\partial_t),d\phi_t),d\phi_t\rangle}_{=-\langle A(d\phi_t(\partial_t),d\phi_t),d\phi_t\rangle=0} \big)\dv\big|_{t=0}\\
=&\int_M\langle\nabla^{\scriptscriptstyle LC}_{e_i}d\phi_t(\partial_t),d\phi_t(e_i)\rangle\dv\big|_{t=0}\\
=&-\int_M\langle\eta,\tau(\phi)\rangle\dv
\end{align*}
confirming the fact that the critical points of the Dirichlet energy are standard harmonic maps even
if we consider a connection with torsion on the target.

Let us now analyze the structure of the torsion endomorphism in \eqref{harmonic-torsion}
in more detail.
Using the decomposition of the torsion endomorphism given in \eqref{torsion-decomposition}
we find
\begin{align*}
A(d\phi,d\phi)=|d\phi|^2V-\langle V,d\phi\rangle d\phi
+T(d\phi,d\phi,\cdot)^\sharp+S(d\phi,d\phi,\cdot)^\sharp.
\end{align*}
Note that due to symmetry reasons the totally antisymmetric part of the torsion does not contribute
to \eqref{harmonic-torsion} as we have
\begin{align*}
T(d\phi,d\phi,\cdot)^\sharp=&T(d\phi(e_i),d\phi(e_j),\cdot)^\sharp g^{ij} \\
=&-T(d\phi(e_j),d\phi(e_i),\cdot)^\sharp g^{ji} \\
=&-T(d\phi,d\phi,\cdot)^\sharp \\
=&0.
\end{align*}
Hence, in general, the equation for a harmonic map with torsion acquires the form
\begin{align*}
\tau(\phi)+|d\phi|^2 V-\langle V,d\phi\rangle d\phi
+S(d\phi,d\phi,\cdot)^\sharp=0.
\end{align*}
In the case of a two-dimensional target only the vectorial torsion contributes
in \eqref{harmonic-torsion} and the equation for harmonic maps with torsion 
further simplifies to
\begin{align*}
\tau(\phi)+|d\phi|^2 V-\langle V,d\phi\rangle d\phi=0.
\end{align*}

Choosing local geodesic normal coordinates \((U,x^i)\) on \(M\) and 
local coordinates \((V,y^\alpha)\) on \(N\) such that \(\phi(U)\subset V\)
the equation for harmonic maps with torsion acquires the form
\begin{align}
\label{euler-lagrange-local}
\Delta \phi^\alpha
+\Gamma^\alpha_{\beta\gamma}(\phi)\frac{\partial\phi^\beta}{\partial x^i}
\frac{\partial\phi^\gamma}{\partial x^j}g^{ij}
+A^\alpha_{\beta\gamma}(\phi)\frac{\partial\phi^\beta}{\partial x^i}
\frac{\partial\phi^\gamma}{\partial x^j}g^{ij}
=0,
\end{align}
where \(\alpha=1,\ldots,\dim N\).

For most of the computations carried out in this article the precise
structure of the torsion endomorphism \(A(X,Y)\) will not be important
and we will mainly work with equation \eqref{harmonic-torsion}.

\subsection{Geodesics with torsion}
As in the case of standard geodesics it is straightforward to see 
that geodesics with torsion have constant speed.
Let \(\gamma\colon I\subset \R\to N\) be a curve
and denote the derivative with respect to the curve parameter by \('\).

\begin{Lem}
Let \(\gamma\colon I\to N\) be a solution of \eqref{harmonic-torsion}.
Then \(|\gamma'|^2\) is constant.
\end{Lem}
\begin{proof}
We compute
\begin{align*}
\frac{d}{ds}|\gamma'|^2=&2\langle\frac{\nabla^{\scriptscriptstyle LC}}{ds}\gamma',\gamma'\rangle
=2\langle\tau(\gamma),\gamma'\rangle
=-2\langle A(\gamma',\gamma'),\gamma'\rangle=0,
\end{align*}
where we first used \eqref{harmonic-torsion} and afterwards the skew-adjointness of the torsion endomorphism 
completing the proof.
\end{proof}
This statement also holds true for magnetic geodesics \cite{MR3778117}
which have some similarity with geodesics with torsion.
In particular, magnetic geodesics also do not always arise from a variational principle
leading to the same technical difficulties as outlined in this article.

In terms of coordinates the equation for geodesics with torsion is given by
\begin{align*}
\gamma''^\alpha=-\big(\Gamma^\alpha_{\beta\delta}(\gamma)+A^\alpha_{\beta\delta}(\gamma)\big)\gamma'^\beta\gamma'^\delta.
\end{align*}
If we interpret this equation as an ordinary differential equation then we know that we can always
extend a solution beyond a given interval of existence as the right hand side is bounded.
However, we cannot expect to find a generalization of the Hopf-Rinow theorem
from Riemannian geometry as was demonstrated in \cite[Section 4]{MR2103403}.

\section{Geometric aspects of harmonic maps with torsion}
In this section we study various geometric aspects related to harmonic maps with torsion.
\subsection{Bochner formulas}
The Bochner technique is a fundamental tool in the analysis of geometric partial differential equations.
It is based on interchanging covariant derivatives and making use of the resulting
curvature terms in order to get a deeper understanding of the solution space of the corresponding
geometric partial differential equation.

In the case of harmonic maps with torsion we can make use of the Bochner technique based on either the Levi-Civita connection or a metric connection with torsion.

Throughout this section
we choose a local orthonormal basis around a point \(p\in M\) such that \(\nabla_{e_i}e_j=0\)
at \(p\) for \(i,j=1,\ldots,\dim M\).

Let us first recall a possible derivation of the Bochner-formula for standard harmonic maps
making use of the Levi-Civita connection.
\begin{align}
\label{bochner-levi-civita}
\Delta^{\scriptscriptstyle LC}d\phi(e_j)
=&\nabla^{\scriptscriptstyle LC}_{e_i}\nabla^{\scriptscriptstyle LC}_{e_i}d\phi(e_j)\\
\nonumber=&
\nabla^{\scriptscriptstyle LC}_{e_i}\nabla^{\scriptscriptstyle LC}_{e_j}d\phi(e_i)\\
\nonumber=&
R^{T^\ast M\otimes\phi^\ast TN}(e_i,e_j)d\phi(e_i)
+\nabla^{\scriptscriptstyle LC}_{e_j}\nabla^{\scriptscriptstyle LC}_{e_i}d\phi(e_i)\\
\nonumber=&
d\phi(\operatorname{Ric}^M(e_j))+R^N(d\phi(e_i),d\phi(e_j))d\phi(e_i)
+\nabla^{\scriptscriptstyle LC}_{e_j}\tau(\phi). 
\end{align}
Note that after the second equals sign we made use of the fact that we employ the Levi-Civita connection.
Moreover, as \(d\phi\in\Gamma(T^\ast M\otimes\phi^\ast TN)\) interchanging covariant derivatives produces the above
curvature term.

Making use of the Bochner formula for the Levi-Civita connection \eqref{bochner-levi-civita} we obtain
\begin{align}
\label{bochner-dphi-lc}
\Delta\frac{1}{2}|d\phi|^2=&|\nabla^{\scriptscriptstyle LC} d\phi|^2+\langle d\phi(\operatorname{Ric}^M),d\phi\rangle
-\langle R^N(d\phi(e_i),d\phi(e_j))d\phi(e_j),d\phi(e_i)\rangle \\
\nonumber&+\langle \nabla^{\scriptscriptstyle LC}_{e_j}\tau(\phi),d\phi(e_j)\rangle.
\end{align}
In the case that \(M\) is compact one can deduce with the help of the maximum principle
that if \(M\) has positive Ricci curvature and \(N\)
non-positive sectional curvature then every harmonic map must be constant.
In the following we want to analyze if this statement remains true in the presence
of torsion on \(N\).

To this end, we derive a generalization of the Bochner-formula \eqref{bochner-levi-civita} taking into account the non-vanishing torsion.

\begin{Prop}[Bochner formula with torsion]
Let \(\phi\colon M\to N\) be a smooth map and suppose that \(N\)
is equipped with a connection with non-vanishing torsion.
Then the following Bochner formula holds
\begin{align}
\label{bochner-torsion}
\Delta^{\scriptscriptstyle Tor}d\phi(e_j)=&
d\phi(\operatorname{Ric}^M(e_j))+R^N_{\scriptscriptstyle Tor}(d\phi(e_i),d\phi(e_j))d\phi(e_i)
+\nabla^{\scriptscriptstyle Tor}_{e_j}\tau^{\scriptscriptstyle Tor}(\phi) \\
\nonumber&+(\nabla^{\scriptscriptstyle LC}_{d\phi(e_i)}A)(d\phi(e_i),d\phi(e_j))
+A(\tau(\phi),d\phi(e_j))\\
&\nonumber+A(d\phi(e_i),\nabla^{\scriptscriptstyle LC}_{e_i}d\phi(e_j))
+A(d\phi(e_i),A(d\phi(e_i),d\phi(e_j)))\\
\nonumber&-(\nabla^{\scriptscriptstyle LC}_{d\phi(e_i)}A)(d\phi(e_j),d\phi(e_i))
-A(\nabla^{\scriptscriptstyle LC}_{e_i}d\phi(e_j),d\phi(e_i))\\
&\nonumber-A(d\phi(e_j),\tau(\phi))
-A(d\phi(e_i),A(d\phi(e_j),d\phi(e_i))).
\end{align}

\end{Prop}
\begin{proof}
First of all, we compute
\begin{align*}
\Delta^{\scriptscriptstyle Tor}d\phi(e_j)
=&\nabla^{\scriptscriptstyle Tor}_{e_i}\nabla^{\scriptscriptstyle Tor}_{e_i}d\phi(e_j)\\
=&\nabla^{\scriptscriptstyle Tor}_{e_i}\big(T(d\phi(e_i),d\phi(e_j))\big)
+\nabla^{\scriptscriptstyle Tor}_{e_i}\nabla^{\scriptscriptstyle Tor}_{e_j}d\phi(e_i)\\
=&\nabla^{\scriptscriptstyle Tor}_{e_i}\big(A(d\phi(e_i),d\phi(e_j))\big)
-\nabla^{\scriptscriptstyle Tor}_{e_i}\big(A(d\phi(e_j),d\phi(e_i))\big)\\
&+R^{T^\ast M\otimes\phi^\ast TN}(e_i,e_j)d\phi(e_i)
+\nabla^{\scriptscriptstyle Tor}_{e_j}\nabla^{\scriptscriptstyle Tor}_{e_i}d\phi(e_i)\\
=&\nabla^{\scriptscriptstyle Tor}_{e_i}\big(A(d\phi(e_i),d\phi(e_j))\big)
-\nabla^{\scriptscriptstyle Tor}_{e_i}\big(A(d\phi(e_j),d\phi(e_i))\big)\\
&+d\phi(\operatorname{Ric}^M(e_j))+R^N_{\scriptscriptstyle Tor}(d\phi(e_i),d\phi(e_j))d\phi(e_i)
+\nabla^{\scriptscriptstyle Tor}_{e_j}\tau^{\scriptscriptstyle Tor}(\phi).
\end{align*}
After the second equals sign we used that we have a connection with metric torsion on 
\(\phi^\ast TN\) and inserted the definition of the torsion tensor \(T(X,Y)\), that is \eqref{def-torsion-tensor}, afterwards.
The result now follows from calculating
\begin{align*}
\nabla^{\scriptscriptstyle Tor}_{e_i}\big(A(d\phi(e_i),d\phi(e_j))\big)
=&(\nabla^{\scriptscriptstyle LC}_{d\phi(e_i)}A)(d\phi(e_i),d\phi(e_j))
+A(\tau(\phi),d\phi(e_j))\\
&+A(d\phi(e_i),\nabla^{\scriptscriptstyle LC}_{e_i}d\phi(e_j))
+A(d\phi(e_i),A(d\phi(e_i),d\phi(e_j)))
\end{align*}
and similarly for the remaining term.
\end{proof}

\begin{Lem}
Let \(\phi\colon M\to N\) be a smooth harmonic map with torsion.
Then the following formula holds
\begin{align}
\label{bochner-dphi-torsion}
\Delta\frac{1}{2}|d\phi|^2=&|\nabla^{\scriptscriptstyle LC} d\phi|^2+\langle d\phi(\operatorname{Ric}^M),d\phi\rangle
-\langle R^N(d\phi(e_i),d\phi(e_j))d\phi(e_j),d\phi(e_i)\rangle \\
\nonumber&
-\langle (\nabla^{\scriptscriptstyle LC}_{d\phi(e_j)}A)(d\phi,d\phi),d\phi(e_j)\rangle
-\langle A(d\phi(e_i),\nabla^{\scriptscriptstyle LC}_{e_j}d\phi(e_i)),d\phi(e_j)\rangle.
\end{align}
\end{Lem}

\begin{proof}
Using that \(\phi\) is a smooth solution of \eqref{harmonic-torsion} we calculate
\begin{align*}
\nabla^{\scriptscriptstyle LC}_{e_j}\tau(\phi)&=
-\nabla^{\scriptscriptstyle LC}_{e_j}\big(A(d\phi,d\phi)\big)\\
&=-(\nabla^{\scriptscriptstyle LC}_{d\phi(e_j)}A)(d\phi,d\phi)
-A(\nabla^{\scriptscriptstyle LC}_{e_j}d\phi(e_i),d\phi(e_i))
-A(d\phi(e_i),\nabla^{\scriptscriptstyle LC}_{e_j}d\phi(e_i)).
\end{align*}
Note that
\begin{align*}
\langle A(\nabla^{\scriptscriptstyle LC}_{e_j}d\phi(e_i),d\phi(e_i)),d\phi(e_j)\rangle
=&\langle A(\nabla^{\scriptscriptstyle LC}_{e_i}d\phi(e_j),d\phi(e_i)),d\phi(e_j)\rangle \\
=&-\langle A(\nabla^{\scriptscriptstyle LC}_{e_i}d\phi(e_j),d\phi(e_j)),d\phi(e_i)\rangle \\
=&-\langle A(\nabla^{\scriptscriptstyle LC}_{e_j}d\phi(e_i),d\phi(e_i)),d\phi(e_j)\rangle=0,
\end{align*}
where we first used the symmetry of the Levi-Civita connection and the skew-symmetry of the torsion endomorphism in the second step. 
Inserting into \eqref{bochner-dphi-lc} completes the proof.
\end{proof}

By a direct calculation we find the following relation between the connection Laplacians 
for the Levi-Civita connection \(\Delta^{\scriptscriptstyle LC}\) and for a connection with torsion \(\Delta^{\scriptscriptstyle tor}\) on \(\phi^\ast TN\)

\begin{align}
\label{relation-laplace}
\Delta^{\scriptscriptstyle Tor}d\phi(e_j)=&\Delta^{\scriptscriptstyle LC}d\phi(e_j)
+A(d\phi,\nabla^{\scriptscriptstyle LC}d\phi(e_j))
+(\nabla_{d\phi(e_i)}^{\scriptscriptstyle LC} A)(d\phi(e_i),d\phi(e_j))\\
\nonumber&+A(\tau(\phi),d\phi(e_j))
+A(d\phi,\nabla^{\scriptscriptstyle LC}d\phi(e_j))
+A(d\phi,A(d\phi,d\phi(e_j))).
\end{align}

\begin{Bem}
\begin{enumerate}
 \item By a direct, but lengthy calculation using \eqref{curvature-tensors} and \eqref{relation-laplace}
it can be verified that \eqref{bochner-torsion} reduces to \eqref{bochner-levi-civita}
if one expresses all quantities in terms of the Levi-Civita connection.
\item It seems rather difficult to extract more information from \eqref{bochner-dphi-torsion} due to the presence of the
torsion terms.
\end{enumerate}
\end{Bem}

\subsection{The effects of conformal transformations}
In this subsection we want to understand the effects of conformal transformations
on both domain and target manifold in the context of harmonic maps with torsion.

In order to analyze the effect of a conformal transformation on the domain
let \(\phi\colon(M,g)\to (N,h)\) be a smooth map.
If we perform a conformal transformation of the metric on the domain,
that is \(\tilde g=e^{2u}g\) for some smooth function \(u\),
we have the following formula for the transformation of the tension field
\begin{align*}
\tau_{\tilde g}(\phi)=e^{-2u}\big(\tau_g(\phi)+(m-2)d\phi(\nabla u)\big),
\end{align*}
where \(\tau_{\tilde g}(\phi)\) denotes the tension field of the map \(\phi\)
with respect to the metric \(\tilde g\) and the Levi-Civita connection on \(N\).

Moreover, it is easy to check that the torsion endomorphism satisfies
\begin{align*}
A_{\tilde g}(d\phi,d\phi)=e^{-2u}A_g(d\phi,d\phi),
\end{align*}
where the notation \(A_g(d\phi,d\phi)\) highlights that we are using the metric \(g\).
Hence, we obtain
\begin{align*}
\tau^{\scriptscriptstyle Tor}_{\tilde g}(\phi)=e^{-2u}\big(\tau^{\scriptscriptstyle Tor}_g(\phi)+(m-2)d\phi(\nabla u)\big).
\end{align*}
We may conclude that (as in the case of standard harmonic maps)
the equation for harmonic maps with torsion \eqref{harmonic-torsion}
is invariant under conformal transformations on the domain in the case that \(\dim M=2\).
The invariance under conformal transformations for standard harmonic maps in \(\dim M=2\) gives rise to 
a conserved quantity, namely that the associated stress-energy tensor is tracefree.
However, in order to obtain the stress-energy tensor one needs to have an energy functional at hand
which one can vary with respect to the metric on the domain. 
As harmonic maps with torsion are non-variational
we cannot expect to find a generalization of the stress-energy tensor for them
and also do not obtain a conserved quantity.

As a second step, we point out an interesting similarity between 
harmonic maps with vectorial torsion and conformal transformations on the target manifold.

One could ask the natural question:
\emph{Suppose we have given a standard harmonic map, can we deform it into
a harmonic map with torsion by performing a conformal transformation
of the metric on the target?}

To approach this question 
let \(\tilde h=e^{2v}h\) be a conformal transformation of the metric on the 
target manifold, where \(v\in C^\infty(N)\). 
Then, for all \(X,Y\in\Gamma(TN)\) the Levi-Civita connections of \(h\) and \(\tilde h\) are related via
\begin{align*}
\tilde\nabla^{\scriptscriptstyle LC}_XY
=\nabla^{\scriptscriptstyle LC}_XY
+X(v)Y+Y(v)X-h(X,Y)\nabla v.
\end{align*}

Hence, the tension fields (computed with respect to the Levi-Civita connection)
satisfy
\begin{align}
\label{tension-field-conformal-target}
\tilde\tau^{\scriptscriptstyle LC}(\phi)=\tau^{\scriptscriptstyle LC}(\phi)
+2\langle d\phi,\nabla v\rangle d\phi-|d\phi|^2\nabla v.
\end{align}
Note that the torsion endomorphism \(A(d\phi,d\phi)\) does not depend on the metric
on the target and is not affected by a conformal transformation on the target.

Now, recall that in the case of vectorial torsion the equation for harmonic maps with
torsion \eqref{harmonic-torsion} acquires the form
\begin{align}
\label{harmonic-map-torsion-vectorial-conformal}
\tau(\phi)+|d\phi|^2 V-\langle V,d\phi\rangle d\phi=0.
\end{align}
In the case that \(V\) is a gradient vector field, that is \(V=\nabla v\)
there is an interesting similarity (up to a factor of two) between \eqref{tension-field-conformal-target}
and \eqref{harmonic-map-torsion-vectorial-conformal}.

\subsection{The stability of harmonic maps with torsion}
In order to understand the stability of standard harmonic maps one usually calculates the 
second variation of the Dirichlet energy \eqref{energy} for a map between two Riemannian manifolds and evaluates it at a critical point,
see for example \cite[Section 1.4.3]{MR1391729} for an overview.
A critical point of the Dirichlet energy, which corresponds
to a standard harmonic map, is stable if the second variation is positive.
Another variant of defining the stability of standard harmonic maps
is to study the spectrum of the associated Jacobi operator. This operator is a linear second
order elliptic operator and can be derived by linearizing the standard harmonic map equation.

In order to discuss the stability of harmonic maps with torsion we cannot consider the second variation of an energy functional.
Nevertheless, we can calculate the linearization of the equation for harmonic
maps with torsion and obtain a corresponding Jacobi operator.

\begin{Satz}
Let \(\phi\colon M\to N\) be a smooth harmonic map with torsion.
In terms of the Levi-Civita connection the corresponding Jacobi operator is given by
\begin{align}
\label{jacobi-lc}
J^{\scriptscriptstyle LC}_\phi(\eta):=\Delta^{\scriptscriptstyle LC}\eta
+R^N(\eta,d\phi)d\phi
+(\nabla^{\scriptscriptstyle LC}_\eta A)(d\phi,d\phi)
+A(\nabla^{\scriptscriptstyle LC}\eta,d\phi)
+A(d\phi,\nabla^{\scriptscriptstyle LC}\eta).
\end{align}
Expressed via a metric connection with torsion the Jacobi-field equation acquires the form
\begin{align}
\label{jacobi-torsion}
J^{\scriptscriptstyle Tor}_\phi(\eta):=&\Delta^{\scriptscriptstyle Tor}\eta
+R^N_{\scriptscriptstyle Tor}(\eta,d\phi)d\phi\\
\nonumber&+(\nabla^{\scriptscriptstyle LC}_{d\phi(e_i)}A)(\eta,d\phi(e_i))
+A(\nabla^{\scriptscriptstyle LC}\eta,d\phi)
+A(\eta,\tau(\phi))
+A(d\phi,A(\eta,d\phi))\\
\nonumber&-(\nabla^{\scriptscriptstyle LC}_{d\phi(e_i)}A)(d\phi(e_i),\eta)
-A(\tau(\phi),\eta)
-A(d\phi,\nabla^{\scriptscriptstyle LC}\eta)
-A(d\phi,A(d\phi,\eta)).
\end{align}

\end{Satz}

\begin{proof}
Suppose we have a smooth harmonic map with torsion,
that is a solution of \eqref{harmonic-torsion}.
To derive the associated Jacobi-field equation
we consider a variation of the map \(\phi\) defined by
\(\phi_t\colon (-\epsilon,\epsilon)\times M\to N\) satisfying 
\(d\phi_t(\partial_t)|_{t=0}=\eta\).

In order to obtain the Jacobi-field equation expressed via the Levi-Civita connection
we calculate
\begin{align*}
0=\frac{\nabla}{\partial t}\big(\tau(\phi_t)+A(d\phi_t,d\phi_t)\big)\big|_{t=0}.
\end{align*}
It is well-known that (see for example \cite[Section 1.4.3]{MR1391729})
\begin{align*}
\frac{\nabla}{\partial t}\tau(\phi_t)\big|_{t=0}=\Delta^{\scriptscriptstyle LC}\eta+R^N(\eta,d\phi)d\phi.
\end{align*}
Regarding the torsion endomorphism we find using the Levi-Civita connection
\begin{align*}
\frac{\nabla}{\partial t}\big(A(d\phi_t,d\phi_t)\big)\big|_{t=0}
=(\nabla^{\scriptscriptstyle LC}_\eta A)(d\phi,d\phi)
+A(\nabla^{\scriptscriptstyle LC}\eta,d\phi)+A(d\phi,\nabla^{\scriptscriptstyle LC}\eta),
\end{align*}
which yields the first statement.

As a second step we derive the Jacobi-field equation making use of a metric  connection with torsion.
Thus, we again compute
\begin{align*}
0=&\frac{\nabla}{\partial t}\tau^{\scriptscriptstyle Tor}(\phi)\big|_{t=0}\\
=&R^N_{\scriptscriptstyle Tor}(\eta,d\phi)d\phi+\nabla_{e_i}^{\scriptscriptstyle Tor}\nabla^{\scriptscriptstyle Tor}_{\partial_t}d\phi_t(e_i)\big|_{t=0} \\
=&R^N_{\scriptscriptstyle Tor}(\eta,d\phi)d\phi
+\Delta^{\scriptscriptstyle Tor}\eta+\nabla_{e_i}^{\scriptscriptstyle Tor}\big(T(d\phi_t(\partial_t),d\phi_t(e_i)\big)\big|_{t=0}\\
=&R^N_{\scriptscriptstyle Tor}(\eta,d\phi)d\phi
+\Delta^{\scriptscriptstyle Tor}\eta
+\nabla_{e_i}^{\scriptscriptstyle LC}\big(A(d\phi_t(\partial_t),d\phi_t(e_i)\big)\big|_{t=0}
+A(d\phi,A(\eta,d\phi))\\
&-\nabla_{e_i}^{\scriptscriptstyle LC}\big(A(d\phi(e_i),d\phi_t(\partial_t))\big)\big|_{t=0}
-A(d\phi,A(d\phi,\eta)).
\end{align*}
By a direct calculation we find
\begin{align*}
\nabla_{e_i}^{\scriptscriptstyle LC}\big(A(d\phi_t(\partial_t),d\phi_t(e_i)\big)\big|_{t=0}=&
(\nabla^{\scriptscriptstyle LC}_{d\phi(e_i)}A)(\eta,d\phi(e_i))
+A(\nabla^{\scriptscriptstyle LC}\eta,d\phi)
+A(\eta,\tau(\phi))
\end{align*}
and a similar formula holds for the remaining term. The claim then follows from combining the 
different equations.
\end{proof}

\begin{Bem}
It is not hard to see that equations \eqref{jacobi-lc} and \eqref{jacobi-torsion}
can easily be transformed into each other.
Using the formula for the connection Laplacians \eqref{relation-laplace} we find
\begin{align*}
\Delta^{\scriptscriptstyle Tor}\eta=&\Delta^{\scriptscriptstyle LC}\eta
+A(d\phi,\nabla^{\scriptscriptstyle LC}\eta)
+(\nabla_{d\phi(e_i)}^{\scriptscriptstyle LC} A)(d\phi(e_i),\eta)\\
&+A(\tau(\phi),\eta)
+A(d\phi,\nabla^{\scriptscriptstyle LC}\eta)
+A(d\phi,A(d\phi,\eta)).
\end{align*}
Moreover, from \eqref{curvature-tensors} we obtain
\begin{align*}
R^N_{\scriptscriptstyle Tor}(\eta,d\phi)d\phi=&R^N_{\scriptscriptstyle LC}(\eta,d\phi)d\phi
+(\nabla^{\scriptscriptstyle LC}_\eta A)(d\phi,d\phi)-(\nabla^{\scriptscriptstyle LC}_{d\phi(e_i)}A)(\eta,d\phi(e_i)) \\
 &+A(\eta,A(d\phi,d\phi))-A(d\phi,A(\eta,d\phi)).
\end{align*}
Inserting both identities into \eqref{jacobi-torsion} then leads to \eqref{jacobi-lc}.
\end{Bem}

As in the case of the standard Jacobi field equation both \eqref{jacobi-lc} and \eqref{jacobi-torsion} constitute a linear elliptic differential operator of second order.
Hence, in the case that the manifold \(M\) is compact we know that both \(J^{\scriptscriptstyle LC}_\phi	\) and \(J^{\scriptscriptstyle Tor}_\phi\) have a discrete spectrum.
For this reason we can in principle perform, as in the case of standard harmonic maps, a stability analysis of harmonic maps with torsion employing methods from spectral geometry.
However, as we can define the Jacobi operator for either the Levi-Civita connection or for a connection with torsion,
there does not seem to be a unique way of defining the stability of a harmonic map with torsion.

It is well known that standard harmonic maps from compact domains to target spaces with negative curvature are stable.
From a spectral point of view this is reflected in the fact that the corresponding Jacobi operator
does not have any eigenvalues.
If we carry out the same analysis for harmonic maps with torsion then it is hard to make
a general statement on the spectrum of the Jacobi operator.
Consequently, in general, we cannot conclude that harmonic maps with torsion to target spaces
with negative curvature (both defined for the Levi-Civita connection or for a connection with torsion) are stable.

In summary, we can say that it is very difficult to understand the stability of harmonic
maps with torsion such that this topic deserves a thorough future investigation.
  
\section{Analytic aspects of harmonic maps with torsion}
This section is devoted to various analytic aspects of a given harmonic map
with torsion.

\subsection{The unique continuation property for harmonic maps with torsion}
It was proved by Sampson \cite[Theorem 2]{MR510549} that standard
harmonic maps satisfy the unique continuation property
and we now extend this result to harmonic maps with torsion.
To this end we recall the following \cite[p. 248]{MR92067}
\begin{Satz}[Aronszajn, 1957]
\label{aro-theorem}
Let \(A\) be a linear elliptic second-order differential operator defined on a domain \(D\) of \(\R^m\).
Let \(u=(u^1,\ldots,u^m)\) be functions in \(D\) satisfying the inequality
\begin{align}
\label{aro-voraus}
|Au^\alpha|\leq C\big(\sum_{i,\alpha}\big|\frac{\partial u^\alpha}{\partial x^i}\big|
+\sum_\alpha|u^\alpha|\big).
\end{align}
If \(u=0\) in an open set, then \(u=0\) throughout \(D\).
\end{Satz}

Making use of this result we can prove the following
\begin{Satz}
Let \(\phi,\phi'\in C^2(M,N)\) be two harmonic maps with torsion.
If \(\phi\) and \(\phi'\) are equal on a connected open set \(W\) of \(M\)
then they coincide on the whole connected component of \(M\) which contains \(W\).
\end{Satz}
\begin{proof}
Let \((U,x^i)\) be a coordinate ball on \(M\), 
by shrinking \(U\) if necessary we can assume that both \(\phi\) and \(\phi'\) map \(U\)
into a single coordinate chart \((V,y^\alpha)\) of \(N\). 
We set \(u^\alpha:=\phi^\alpha-\phi'^\alpha\) and 
using the local form of the Euler-Lagrange equation \eqref{euler-lagrange-local} we find
\begin{align*}
\Delta u^\alpha=&-\Gamma^\alpha_{\beta\gamma}(\phi)\frac{\partial \phi^\beta}{\partial x^i}\frac{\partial \phi^\gamma}{\partial x^j}g^{ij}
+\Gamma^\alpha_{\beta\gamma}(\phi')\frac{\partial \phi'^\beta}{\partial x^i}\frac{\partial \phi'^\gamma}{\partial x^j}g^{ij} \\
&-A_{\beta\gamma}^\alpha(\phi)\frac{\partial \phi^\beta}{\partial x^i}\frac{\partial \phi^\gamma}{\partial x^j}g^{ij}
+A_{\beta\gamma}^\alpha(\phi')\frac{\partial \phi'^\beta}{\partial x^i}\frac{\partial \phi'^\gamma}{\partial x^j}g^{ij}\\
=&-\big(\Gamma^\alpha_{\beta\gamma}(\phi)-\Gamma^\alpha_{\beta\gamma}(\phi')\big)
\frac{\partial\phi^\beta}{\partial x^i}\frac{\partial\phi^\gamma}{\partial x^j}g^{ij}
-\Gamma^\alpha_{\beta\gamma}(\phi')\frac{\partial u^\beta}{\partial x^i}\frac{\partial\phi^\gamma}{\partial x^j}g^{ij}
-\Gamma^\alpha_{\beta\gamma}(\phi')\frac{\partial \phi^\beta}{\partial x^i}\frac{\partial u^\gamma}{\partial x^j}g^{ij}\\
&-\big(A^\alpha_{\beta\gamma}(\phi)-A^\alpha_{\beta\gamma}(\phi')\big)
\frac{\partial \phi^\beta}{\partial x^i}\frac{\partial \phi^\gamma}{\partial x^j}g^{ij}
-A^\alpha_{\beta\gamma}(\phi')\frac{\partial u^\beta}{\partial x^i}\frac{\partial\phi^\gamma}{\partial x^j}g^{ij}
-A^\alpha_{\beta\gamma}(\phi')\frac{\partial \phi^\beta}{\partial x^i}\frac{\partial u^\gamma}{\partial x^j}g^{ij}\\
\leq &C(\sum_{i,\alpha}\big|\frac{\partial u^\alpha}{\partial x^i}\big|+\sum_\alpha|u^\alpha|),
\end{align*}
where we used the mean-value inequality to estimate the terms involving the Christoffel symbols and the torsion endomorphisms.
Moreover, we utilized that \(\phi\) and its derivatives
are bounded as we are considering maps between compact sets.
The result now follows by application of Theorem \ref{aro-theorem}. 
\end{proof}

\begin{Bem}
Note that the same proof from above would also work for VT-harmonic maps \eqref{vt-harmonic}
which are a variant of harmonic maps with torsion that also contain a term that is linear in \(d\phi\)
on the right hand side.
\end{Bem}

\subsection{Regularity of weak solutions}
In this subsection we establish the regularity of weak harmonic maps with torsion.
It turns out that the additional nonlinearity arising due to the non-vanishing torsion
has the right antisymmetric structure such that the known regularity theory for standard harmonic maps is still applicable.

First, we define a weak solution of \eqref{harmonic-torsion}.
\begin{Dfn}
We call a map \(\phi\in W^{1,2}(M,N)\) a \emph{weak harmonic map with torsion} if it
solves \eqref{harmonic-torsion} in the sense of distributions.
\end{Dfn}

We will prove the following regularity result for weak harmonic maps with torsion.
\begin{Satz}
\label{theorem-regularity}
Let \((M,g)\) be a closed Riemannian manifold with \(\dim M=m\) and \((N,h)\)
a compact Riemannian manifold.
Suppose that \(\phi\colon M\to N\) is a weak solution of \eqref{harmonic-torsion}
that satisfies the Morrey growth condition
\begin{align}
\label{morrey-growth-theorem}
\sup_{x\in U}\big(\frac{1}{r^{m-2}}\int_{B_r(x)\cap U}|d\phi|^2d\mu\big)^\frac{1}{2}
<\epsilon
\end{align}
for some small \(\epsilon>0\).

Then a weak solution \(\phi\in W^{1,2}(M,N)\) is smooth in \(U\), where \(U\) is an open subset of \(M\).
\end{Satz}

\begin{Bem}
\begin{enumerate}
\item In the case that \(M\) is a closed surface it is enough to demand that \(\phi\in W^{1,2}(M,N)\)
as we can always find a small disc on which \eqref{morrey-growth-theorem} holds.

\item In the regularity analysis of standard harmonic maps one usually demands that the harmonic map
is stationary which is a weaker condition as \eqref{morrey-growth-theorem}.
A harmonic map is stationary if it is a critical point of the energy \(E(\phi)\) both with respect
to variations of the map as well as with respect to variations of the metric on the domain.
The condition of being stationary can be interpreted as a weak formulation of the invariance
of the energy under diffeomorphisms on the domain \cite[Section 2.3]{MR4007262}.
As harmonic maps with torsion do not allow for a variational formulation we cannot come up 
with a generalization of the notion of being stationary
but have to demand the smallness condition \eqref{morrey-growth-theorem}.
\end{enumerate}

\end{Bem}

In order to prove Theorem \ref{theorem-regularity} it is natural to apply the embedding theorem of Nash and to isometrically
embed the target manifold \(N\), where \(n:=\dim N\), into some \(\R^q\) of sufficiently large dimension \(q\).
We denote this isometric embedding by \(\iota\colon N\to\R^q\) which we assume to be smooth.
Now, we consider the composite map \(\phi':=\iota\circ\phi\colon M\to\R^q\) and 
let \(u^\alpha,1\leq\alpha\leq q\) be global coordinates on the ambient space \(\R^q\).
Moreover, let \(\nu_\theta,\theta=n+1,\ldots,q\) be an orthonormal frame of the submanifold \(\iota(\phi)\). In the following we will still write \(\phi\) instead of \(\phi'\)
in order to shorten the notation.

\begin{Lem}
\label{lem-extrinsic-antisymmetric}
The extrinsic version of the equation for harmonic maps with torsion \eqref{harmonic-torsion}
for \(\phi\colon M\to\R^q\) is given by
\begin{align}
\label{extrinsic-phi}
-\Delta\phi^\alpha&=(\omega_{i}^{\alpha\beta}+B_i^{\alpha\beta})\frac{\partial\phi^\beta}{\partial x^i},
\end{align}
where \(1\leq\alpha\leq q\) and
\begin{align*}
\omega_{i}^{\alpha\beta}&:=
\frac{\partial\phi^\gamma}{\partial x^i}\frac{\partial\nu_l^\beta}{\partial u^\gamma}\nu_l^\alpha-\frac{\partial\phi^\gamma}{\partial x^i}\frac{\partial\nu_l^\alpha}{\partial u^\gamma}\nu_l^\beta
=-\omega_{i}^{\alpha\beta},\\
B_i^{\alpha\beta}&:=A_{\gamma\beta}^{\alpha}\frac{\partial\phi^\gamma}{\partial x^i}=-B_i^{\beta\alpha}.
\end{align*}
\end{Lem}
\begin{proof}
It is well known by now that the standard harmonic equation can be written in the form \eqref{extrinsic-phi}, 
see for example the discussion in the introduction of \cite{MR2285745}.
We can extend the torsion endomorphism \(A\) to the ambient space by parallel transport.
Moreover, we set
\begin{align*}
A_{\gamma\beta}^{\alpha}:=h^{\alpha\delta}A_{\gamma\beta\delta}=h^{\alpha\delta}\langle A(\partial_{y^\gamma},\partial_{y^\beta}),\partial_{y^\delta}\rangle,
\end{align*}
where \(\partial_{y^\alpha}\) is a local basis of \(TN\).
The antisymmetry follows from \eqref{torsion-skewadjointess}.
\end{proof}

Now, we recall the following regularity result.
Suppose we have a distributional solution to 
\begin{align}
\label{laplace-antisymmetric}
-\Delta u=\Omega\cdot\nabla u,
\end{align}
where \(u\colon B^m\to\R^q\) and \(\Omega\) being antisymmetric.
Here, \(B^m\) denotes the \(m\)-dimensional unit disc.

For such kind of systems the following regularity result holds \cite[Theorem 1.1]{MR2383929}:
\begin{Satz}
\label{regularity-riviere}
For every \(m\in\N\) there exists \(\epsilon(m)>0\) such that for every
\(\Omega\in L^2(B^m,\mathfrak{so}(q)\otimes\R^m)\)
and for every weak solution \(u\in W^{1,2}(B^m,\R^q)\) of \eqref{laplace-antisymmetric}
that satisfies the Morrey growth condition
\begin{align}
\sup_{x\in B}\big(\frac{1}{r^{m-2}}\int_{B_r(x)\cap B}(|\nabla u|^2+|\Omega|^2)d\mu\big)^\frac{1}{2}<\epsilon(m)
\end{align}
we have that \(u\) is locally Hölder-continuous in \(B^m\)
with exponent \(0<\alpha=\alpha(m)<1\).
\end{Satz}

We are now able to give the proof of Theorem \ref{theorem-regularity}.

\begin{proof}[Proof of Theorem \ref{theorem-regularity}]
Thanks to \eqref{extrinsic-phi} the equation for harmonic maps with
torsion can be written in the form \eqref{laplace-antisymmetric}.
In addition, as \(N\) is compact by assumption and since \(\phi\) is a weak solution of \eqref{harmonic-torsion}
we have that \(\Omega\in L^2(B^m,\mathfrak{so}(q)\otimes\R^m)\).
As also the smallness condition \eqref{morrey-growth-theorem} holds by assumption 
we can apply Theorem \ref{regularity-riviere} to 
the extrinsic version of the equation for harmonic maps with torsion \eqref{extrinsic-phi}
yielding the local Hölder continuity of \(\phi\).
From this point on the smoothness of \(\phi\) follows from the same bootstrap argument as for standard harmonic maps,
see for example \cite[Theorem 7.4.1]{MR2829653}.
\end{proof}

\subsection{Local energy estimates and a removable singularity theorem}
A classic result of Sacks and Uhlenbeck shows that standard harmonic maps of finite Dirichlet energy
from a two-dimensional domain cannot have isolated point singularities or
phrased differently isolated point singularities can always be removed \cite{MR604040}.
One should expect that such a result also holds for harmonic maps with torsion
as their nonlinear structure on the right hand side 
is the same as in the case of standard harmonic maps.

However, it turns out that the method of proof applied by Sacks and Uhlenbeck does not
directly carry over to harmonic maps with torsion. At the heart of their proof
is a Pohozaev identity which is derived from the stress-energy tensor.
As harmonic maps with torsion do not arise from a variational 
principle it is not possible to extend the methods used in the analysis
of standard harmonic maps and we have to make use of a different method of prove here.

In the following it will again turn out to be useful to 
apply the following extrinsic version of the equation for harmonic maps with torsion
\begin{align}
\label{extrinsic-phi-b}
-\Delta\phi=\sff(d\phi,d\phi)+A(d\phi,d\phi).
\end{align}
Here, \(\phi\colon M\to\R^q\) is a vector-valued map and \(\sff\) denotes the second fundamental form
of the embedding into the ambient space \(\R^q\).
To obtain this extrinsic version of the equation for harmonic maps with torsion one again applies
the embedding theorem of Nash, see the discussion in front of Lemma \ref{lem-extrinsic-antisymmetric}
and the introduction of \cite{MR2285745} for more details.

Note that harmonic maps with torsion have the same analytic structure as standard harmonic maps
in the sense that they satisfy
\begin{align*}
\Delta\phi\leq C|d\phi|^2
\end{align*}
for a positive constant \(C>0\).

This suggests to define the following local energy
\begin{Dfn}
The local energy of a map \(\phi\colon M\to N\) is defined as follows
\begin{align}
E(\phi,U)=\int_U|d\phi|^2 d\mu,
\end{align}
where \(U\subset M\).
\end{Dfn}

We obtain the following result

\begin{Satz}[Removable Singularity Theorem]
\label{theorem-removable}
Let \((M,g)\) be a closed Riemannian manifold of dimension \(\dim M=2\)
and \((N,h)\) a compact Riemannian manifold.
Let \(\phi\) be a harmonic map with torsion which is smooth on 
\(U\setminus\{p\}\) for some \(p\in U\subset M\). If \(\phi\)
has finite energy \(E(\phi,U)\), then \(\phi\) extends to a smooth solution on \(U\).
\end{Satz}

In order to prove Theorem \ref{theorem-removable} we set \(D'=D\setminus\{0\}\),
where \(D\subset M\) denotes the unit disc in two dimensions,
and make use of the following

\begin{Lem}
\label{lemma-extension}
Let \((M,g)\) be a closed Riemannian manifold of dimension \(\dim M=2\)
and suppose \(\phi\in W^{1,2}(D',\R^q)\) satisfies
\begin{align}
\label{lemma-weak-assumption}
\int_{D'}\langle\nabla\phi,\nabla\eta\rangle ~d\mu=\int_{D'}f(x,\phi,\nabla\phi)\eta ~d\mu
\end{align}
for all \(\eta\in W^{1,2}_0\cap L^\infty(D',\R^q)\), where \(f\) satisfies the growth assumption
\begin{align*}
|f(x,\phi,p)|\leq a+b|p|^2
\end{align*}
with constants \(a,b>0\) for all \((x,\phi,p)\in D'\times\R^q\times\R^{2q}\).
Then we also have
\begin{align*}
\int_{D}\langle\nabla\phi,\nabla\eta\rangle ~d\mu=\int_{D}f(x,\phi,\nabla\phi)\eta ~d\mu.
\end{align*}
\end{Lem}
\begin{proof}
For a proof see \cite[Lemma A.2, p. 225]{MR1100926}.
\end{proof}

We can now give the proof of Theorem \ref{theorem-removable}.

\begin{proof}[Proof of Theorem \ref{theorem-removable}]
First of all, we note that \eqref{lemma-weak-assumption} is the weak version 
of \eqref{extrinsic-phi-b}. Hence, due to the assumption \(E(\phi,U)<\infty\) we 
can apply Lemma \ref{lemma-extension}.
Thus, we know that every weak harmonic map with torsion
defined on a disc with the origin removed can be extended to a weak harmonic map
with torsion on the whole unit disc. 
However, the assumptions of Theorem \ref{theorem-removable}
allow to also apply the regularity result Theorem \ref{theorem-regularity} such that
all weak solutions considered here are actually smooth.
\end{proof}

Besides the removable singularity theorem we can also give the following
local energy estimates.

\begin{Satz}[\(\epsilon\)-Regularity Theorem]
\label{theorem-epsilon-regularity}
Assume that \(\phi\) is a smooth harmonic map with torsion
with small energy 
\begin{align*}
\label{small-energy}
E(\phi,D)<\epsilon.
\end{align*}
Then the following estimate holds
\begin{align}
\|d\phi\|_{W^{1,p}(\tilde{D})}\leq C(\tilde{D},p)\|d\phi\|_{L^2(D)}
\end{align}
for all \(\tilde{D}\subset D, p>1\), where \(C(\tilde{D},p)\) is a positive
constant depending only on \(\tilde{D}\) and \(p\).
\end{Satz}
\begin{proof}
This can be proved by the same method as for standard harmonic maps,
see for example \cite[Lemma 2.4.1]{MR1100926} or \cite[Theorem 3.3]{MR3558358}.
\end{proof}

Exploiting the scaling invariance of \eqref{extrinsic-phi-b} we also obtain
\begin{Cor}
\label{corollary-energy-local}
There is an \(\epsilon>0\) small enough such that if \(\phi\) is a smooth harmonic map with torsion
with finite energy \(E(\phi,D)<\epsilon\)
then for any \(x\in D_{\frac{1}{2}}\) we have
\begin{align}
|d\phi(x)||x|\leq C\|d\phi\|_{L^2(D_{2|x|})}, 
\end{align}
where \(C\) is a positive constant.
\end{Cor}
\begin{proof}
This follows from a scaling argument, fix any \(x_0\in D\setminus\{0\}\) and define \(\tilde{\phi}\) as
\[
\tilde{\phi}(x):=\phi(x_0+|x_0|x).
\]
It is easy to see that \(\tilde{\phi}\) is a smooth solution of \eqref{extrinsic-phi-b} on \(D\)
with \(E(\tilde{\phi},D)<\epsilon\). By application of Theorem \ref{theorem-epsilon-regularity}, we have
\begin{align*}
|d\tilde{\phi}|_{L^\infty(D_\frac{1}{2})}\leq C\|d\tilde{\phi}\|_{L^2(D)}
\end{align*}
and scaling back yields the assertion.
\end{proof}

\subsection{A Liouville theorem on complete non-compact manifolds}
In this subsection we will establish a vanishing result for solutions of \eqref{harmonic-torsion}
on a large class of complete non-compact domain manifolds.
For standard harmonic maps a corresponding result was obtained in \cite{MR1333821}, see
also \cite{MR3886921} for further generalizations.

At the heart of the proof is a \emph{Euclidean type Sobolev inequality}, that is
\begin{align}
\label{sobolev-inequality}
(\int_M|u|^{2m/(m-2)}\dv)^\frac{m-2}{m}\leq C_2\int_M|\nabla u|^2\dv
\end{align}
for all \(u\in W^{1,2}(M)\) with compact support,
where \(C_2\) is a positive constant that depends on the geometry of \(M\).
Such an inequality holds in \(\R^m\) and is well-known as \emph{Gagliardo-Nirenberg inequality}.
However, if one considers a non-compact complete Riemannian manifold of infinite volume
one has to make additional assumptions to have an equality 
of the form \eqref{sobolev-inequality} at hand, see the introduction of \cite{MR3886921}
for more details.

We will make use of a cutoff function  \(0\leq\eta\leq 1\) on \(M\) that satisfies
\begin{align*}
\eta(x)=1\textrm{ for } x\in B_R(x_0),\qquad \eta(x)=0\textrm{ for } x\in B_{2R}(x_0),\qquad |\nabla\eta|\leq\frac{C}{R}\textrm{ for } x\in M,
\end{align*}
where \(B_R(x_0)\) denotes the geodesic ball around \(x_0\) with radius \(R\).

Finally, we obtain the following Liouville theorem:

\begin{Satz}
\label{theorem-liouville}
Let \((M,g)\) be a complete and non-compact Riemannian manifold of dimension \(\dim M=m>2\)
with positive Ricci curvature that admits a Euclidean type Sobolev inequality of the form \eqref{sobolev-inequality}. 
Moreover, let \(N\) be a Riemannian manifold of bounded geometry, that is 
\(|R^N|_{L^\infty}+|A|_{L^\infty}<\infty\).
Assume that \(\phi\) is a smooth solution of \eqref{harmonic-torsion}.
If
\begin{align}
\label{assumption-smallness}
\int_M|d\phi|^{m}\dv<\epsilon,
\end{align}
where \(\epsilon>0\) satisfies
\begin{align*}
\epsilon\leq \big[\frac{2m-4}{m^2C_s}
\big(|R^N|_{L^\infty}+|A|_{L^\infty}\frac{1+m}{2}\big)^{-1}
\big]^\frac{m}{2}
\end{align*}
then \(\phi\) must be trivial.
\end{Satz}

First, we will derive the following inequality similar to \cite[Lemma 2.2]{MR3886921}.

\begin{Lem}
\label{lem-phi-lr}
Suppose that \(\phi\colon M\to N\) is a smooth solution of \eqref{harmonic-torsion} and 
\(\dim M>2\).
Then the following inequality holds
\begin{align}
\label{ineq-phi}
\nonumber\frac{C}{R^2}\int_M|d\phi|^m\dv\geq&\int_M\eta^2\langle d\phi(\operatorname{Ric}^M),d\phi\rangle|d\phi|^{m-2}\dv 
+\frac{1}{2}\int_M\eta^2|\nabla d\phi|^2|d\phi|^{m-2}\dv \\
&+\frac{2m-4}{m^2}\int_M\eta^2\big|d|d\phi|^\frac{m}{2}\big|^2\dv \\
\nonumber&-\big(|R^N|_{L^\infty}+\frac{1+m}{2}|A|_{L^\infty}\big)\int_M\eta^2|d\phi|^{m+2}\dv.
\end{align}
\end{Lem}
\begin{proof}
Testing the Bochner formula for the Levi-Civita connection \eqref{bochner-levi-civita}
with \(\eta^2|d\phi|^{m-2}d\phi(e_i)\) and integrating over \(M\) we find
\begin{align*}
\int_M\eta^2\langle\Delta d\phi,d\phi\rangle|d\phi|^{m-2}\dv
=&\int_M\eta^2\langle d\phi(\text{Ric}^M(e_i)),d\phi(e_i)\rangle|d\phi|^{m-2}\dv \\
&+\int_M\eta^2\langle R^N(d\phi(e_i),d\phi(e_j))d\phi(e_i),d\phi(e_j)\rangle|d\phi|^{m-2}\dv\\
&+\int_M\eta^2\langle\nabla\tau(\phi),d\phi\rangle|d\phi|^{m-2}\dv.
\end{align*}

Using integration by parts we may rewrite
\begin{align*}
\int_M\eta^2\langle\nabla\tau(\phi),d\phi\rangle|d\phi|^{m-2}\dv=&-2\int_M\eta\nabla\eta\langle\tau(\phi),d\phi\rangle|d\phi|^{m-2}\dv\\
&-\int_M\eta^2|\tau(\phi)|^2|d\phi|^{m-2}\dv \\
&-(m-2)\int_M\eta^2\langle\tau(\phi),d\phi\rangle|d\phi|^{m-4}\langle d\phi,\nabla d\phi\rangle\dv
\end{align*}
and also
\begin{align*}
\int_M\eta^2\langle\Delta d\phi,d\phi\rangle|d\phi|^{m-2}\dv=&-2\int_M\eta\nabla\eta\langle\nabla d\phi,d\phi\rangle|d\phi|^{m-2}\dv 
-\int_M\eta^2|\nabla d\phi|^2|d\phi|^{m-2}\dv \\
&-(m-2)\int_M\eta^2|\langle\nabla d\phi,d\phi\rangle|^2|d\phi|^{m-4}\dv.
\end{align*}
This allows us to deduce the following inequality
\begin{align*}
\int_M\eta^2&\langle d\phi(\text{Ric}^M(e_i)),d\phi(e_i)\rangle|d\phi|^{m-2}\dv
+\int_M\eta^2|\nabla d\phi|^2|d\phi|^{m-2}\dv \\
&+(m-2)\int_M\eta^2|\langle\nabla d\phi,d\phi\rangle|^2|d\phi|^{m-4}\dv \\
\leq &\frac{C}{R^2}\int_M|d\phi|^m\dv+\frac{1}{2}\int_M\eta^2|\nabla d\phi|^2|d\phi|^{m-2}\dv 	
+\frac{1+m}{2}\int_M\eta^2|\tau(\phi)|^2|d\phi|^{m-2}\dv \\
&+|R^N|_{L^\infty}\int_M\eta^2|d\phi|^{m+2}\dv
+\frac{m-2}{2}\int_M\eta^2|\langle\nabla d\phi,d\phi\rangle|^2|d\phi|^{m-4}\dv.
\end{align*}
The preceding inequality yields
\begin{align*}
\int_M\eta^2&\langle d\phi(\text{Ric}^M(e_i)),d\phi(e_i)\rangle|d\phi|^{m-2}\dv
+\frac{1}{2}\int_M\eta^2|\nabla d\phi|^2|d\phi|^{m-2}\dv \\
&+\frac{m-2}{2}\int_M\eta^2|\langle\nabla d\phi,d\phi\rangle|^2|d\phi|^{m-4}\dv \\
\leq &\frac{C}{R^2}\int_M|d\phi|^m\dv
+\frac{1+m}{2} 
+\big(|R^N|_{L^\infty}+\frac{1+m}{2}|A|_{L^\infty}\big)\int_M\eta^2|d\phi|^{m+2}\dv,
\end{align*}
where we have used that \(\phi\) is a solution of \eqref{harmonic-torsion} in the last step.
The claim follows from using
\begin{align*}
|\langle\nabla d\phi,d\phi\rangle|^2|d\phi|^{m-4}
=&\frac{1}{4}\big|d|d\phi|^2\big|^2|d\phi|^{m-4}
=\frac{4}{m^2}\big|d|d\phi|^\frac{m}{2}\big|^2.
\end{align*}
\end{proof}

To complete the proof we also need the following
\begin{Lem}
\label{lem-sobolev}
Let \((M,g)\) be a complete non-compact Riemannian manifold of infinite volume that admits 
a Euclidean type Sobolev inequality.
Assume that \(f\) is a positive function on \(M\).
For \(\dim M=m>2\) the following inequality holds
\begin{align}
\label{inequality-sobolev}
\int_M\eta^2f^{m+2}\dv\leq C_s\big(\int_Mf^m\dv\big)^\frac{2}{m}\big(\frac{1}{R^2}\int_Mf^{m}\dv
+\int_M\eta^2|df^\frac{m}{2}|^2\dv\big),
\end{align}
where the positive constant \(C_s\) depends on \(m\) and the geometry of \(M\).
\end{Lem}
\begin{proof}
For a derivation see the proof of \cite[Lemma 2.4]{MR3886921}.
\end{proof}

At this point we are ready to complete the proof of Theorem \ref{theorem-liouville}.

\begin{proof}[Proof of Theorem \ref{theorem-liouville}]
We apply \eqref{inequality-sobolev} to the last term on the right hand side of \eqref{ineq-phi}
with \(f=|d\phi|\)
and by choosing \(\epsilon\) small enough 
and taking the limit \(R\to\infty\) we obtain from \eqref{ineq-phi} that
\begin{align*}
0\geq&\int_M\langle d\phi(\operatorname{Ric}^M),d\phi\rangle|d\phi|^{m-2}\dv.
\end{align*}
By assumption \(M\) has positive Ricci curvature and thus \(\text{Ric}^M\) is a positive definite non-degenerate bilinear form on \(TM\).
As \(m>2\) this completes the proof.
\end{proof}

\textbf{Acknowledgements:}
The author gratefully acknowledges the support of the Austrian Science Fund (FWF) 
through the project P30749-N35 ``Geometric variational problems from string theory''.

\bibliographystyle{plain}
\bibliography{mybib}
\end{document}